\documentclass[twoside]{article}
\usepackage[a4paper,inner=3cm,outer=3cm]{geometry}
\pdfoutput=1

\usepackage{setspace}
\usepackage{Style}
\usepackage{pdfpages}

  \usepackage[numbers,sort&compress]{natbib}%
\renewcommand{\hc}[1]{\left.\left(#1\right)\right|_{\textnormal{h+c}}}
\newcommand{\hcS}[1]{\left.#1\right|_{\textnormal{h+c}}}
\renewcommand{\hcG}[1]{\left.\left(#1\right)\right|_{\textnormal{h+c}}^{\gamma}}
\renewcommand{\hcZ}[1]{\left.\left(#1\right)\right|_{\textnormal{h+c}}^{\gamma=0}}
\renewcommand{\hcN}[1]{\left.\left(#1\right)\right|_{\textnormal{h+c}}^{\gamma\ne0}}

\title{Partial symmetries of iterated plethysms}
\author{Álvaro Gutiérrez \and Mercedes H. Rosas}

\usepackage{fancyhdr}

\fancypagestyle{alim}{\fancyhf{}\fancyfoot[R]{\footnotesize\vspace{-5.69em}
MR was partially supported by MTM2016-75024-P and FEDER,\\ PID2020-117843GB-I00, and Proyectos I+D+i FEDER Andalucía US-1262169.\\
AG was partially supported by
Junta de Andalucía FQM-333.}}
\begin{document}

\thispagestyle{alim}
\begin{center}\large {\ }
\vspace{1em}
\textbf{\scshape\Large Partial symmetries of iterated plethysms}
\\ \vspace{2em}
\renewcommand{\thefootnote}{\fnsymbol{footnote}}
\renewcommand{\arraystretch}{1}
\setlength{\tabcolsep}{2em}
\begin{tabular}{cc}
    \'Alvaro Gutiérrez\footnotemark[1]{}
    & Mercedes H. Rosas\footnotemark[2]{} \\
    {\small Department of Mathematics} & {\small Departamento de Álgebra}\\
    {\small University of Bonn, Germany} & {\small Universidad de Sevilla, Spain}
\end{tabular}
\footnotetext[1]{gut@uni-bonn.de}
\footnotetext[2]{mrosas@us.es}
\renewcommand{\thefootnote}{\arabic{footnote}}
\setcounter{footnote}{0}
\\ \vspace{2em}
\end{center}

\begin{abstract}This work highlights the existence of  partial symmetries  in large families of iterated plethystic coefficients. The plethystic coefficients involved  come from the  expansion in the Schur basis of   iterated plethysms of Schur functions indexed by  one-row partitions.

The partial symmetries are described in terms of  an involution on partitions, the flip involution, that generalizes  the ubiquitous $\omega$ involution. Schur-positive symmetric functions possessing this partial symmetry are termed flip-symmetric.

The operation of taking plethysm with $s_\lambda$ preserves  flip-symmetry, provided that $\lambda$ is a partition of two.
Explicit formulas for the iterated plethysms $s_2\circ s_b\circ s_a$ and $s_c\circ s_2\circ s_a$, with $a,$  $b,$ and $c$ $\ge$ $2$  allow us to show that these two families of iterated plethysms are flip-symmetric.
The article concludes with some observations, remarks, and open questions on the  unimodality and asymptotic normality of certain
flip-symmetric sequences of iterated plethystic coefficients.\\

\end{abstract}

\tableofcontents

\newpage
\section{Introduction}
\newcommand{\GL}{\textnormal{GL}}
\newcommand{\Gl}{\textnormal{GL}}
Let $V$ be an $n$ dimensional complex vector space.
Each partition $\lambda$ of length at most $n$ indexes an 
irreducible representation (unique up to isomorphism)  of the complex general lineal group $\GL(V).$ 
The irreducible representation   indexed by $\lambda$
can be  constructed as the evaluation of the Schur functor $\SS^\lambda$ on the vector space $V.$ Therefore, it is denoted   by $\SS^\lambda[V]$.

The composition of representations provides us with a important and natural way of combining group representations, an operation  referred to as the plethysm of representations. In the
setting of the representation theory of the general lineal group, the plethysm of the irreducible representations indexed by $\mu$  and $\nu$ is  defined by the composition of Schur functors
$ \SS^\mu[ \SS^\nu[V]]$. 
 Further information  can be found in Fulton and Harris' book \cite{Fulton_Harris}.

Rational representations 
of the general linear group are completely reducible. This raises
 the question of decomposing the plethysm $\SS^\mu[ \SS^\nu[V]]$  as a sum of irreducible representations.
 
The plethystic coefficient $a_{\mu[\nu]}^\lambda$ is defined as the multiplicity of $ \SS^\lambda[V]$ in  $ \SS^\mu[ \SS^\nu[V]]$.
More generally, the iterated plethystic coefficient $a_{\mu^1[\mu^2[\ldots[\mu^k]]]}^\lambda$ is defined as the
multiplicity of $ \SS^\lambda[V]$ in  $\SS^{\mu^1}[\SS^{\mu^2}[\ldots[\SS^{\mu^k}[V]]]]$.  The partitions indexing the iterated plethystic coefficient
$a_{\mu^1[\mu^2[\ldots[\mu^k]]]}^\lambda$ satisfy that $|\lambda| = |\mu^1| |\mu^2| \cdots |\mu^k|$.

The problem of understanding  the plethystic coefficients is a notoriously hard problem \cite{Fulton_Harris, colmenarejo17, fischer, kahle} that has stumped many attempts to solve it.
In this work, we report the occurrence of partial symmetries in certain iterated plethystic coefficients indexed by a specific type of partitions, hook+column partitions. First considered by 
 Langley and Remmel \cite{langley} in 2004, a hook+column
 partition is a partition of the form $(\alpha, 2^\beta, 1^\gamma)$.
Langley and Remmel  obtained  a simple formula (stated in Theorem \ref{sasb}) for the plethystic coefficients
$a^\lambda_{(b)[(a)]}$, where $\lambda$ is a  hook+column partition. 

Following them, we restrict our attention to iterated plethystic coefficients $a_{\mu^1[\mu^2[\ldots[\mu^k]]]}^\lambda$ 
where  $\lambda$ is a  hook+column, and where each partition $\mu^i$ is either a row or a column partition. We derive closed formulas for the 
plethystic coefficients 
$a_{(c)[(b)[(a)]]}^\lambda$  when either $b$ or $c$ is equal to 2, and $\lambda$ is a hook+column partition. These formulas (described in Theorems \ref{s2 o sb o sa} and \ref{sc o s2 o sa}) allow us to uncover new and intriguing partial symmetries in the iterated plethystic coefficients that we illustrate in the following example.

\begin{ex}\label{ex:intro}\yindex
Consider  the coefficients $a_{\ydiagram{2}[\ydiagram{3}[\ydiagram{2}]]}^{\lambda}$ with
$\lambda$ hook+column partitions.  The nonzero coefficients are:

\begin{center}
    {
    \renewcommand{\arraystretch}{2.5}
    \setlength{\tabcolsep}{0.25em}
    \scriptsize 
    \vspace{.5em } 
    \begin{tabular}{cccccc} $a_{\ydiagram{2}[\ydiagram{3}[\ydiagram{2}]]}^{\ydiagram{2, 2, 2, 2, 2, 2}} = 1$, &
        $a_{\ydiagram{2}[\ydiagram{3}[\ydiagram{2}]]}^{\ydiagram{2, 2, 2, 2, 4}} = 2$, &
        $a_{\ydiagram{2}[\ydiagram{3}[\ydiagram{2}]]}^{\ydiagram{2, 2, 2, 6}} = 3$, &
        $a_{\ydiagram{2}[\ydiagram{3}[\ydiagram{2}]]}^{\ydiagram{2, 2, 8}} = 3$, &
        $a_{\ydiagram{2}[\ydiagram{3}[\ydiagram{2}]]}^{\ydiagram{2, 10}} = 2$, &
        $a_{\ydiagram{2}[\ydiagram{3}[\ydiagram{2}]]}^{\ydiagram{12}} = 1$,\\
        & &
        $a_{\ydiagram{2}[\ydiagram{3}[\ydiagram{2}]]}^{\ydiagram{1, 2, 2, 2, 5}} = 1$, &
        $a_{\ydiagram{2}[\ydiagram{3}[\ydiagram{2}]]}^{\ydiagram{1, 2, 2, 7}} = 1$.& &
    \end{tabular}
    \vspace{.5em } 
    }
\end{center}
The symmetry demonstrated in this example can be characterized by a generalization of the transposition of diagrams where some $(2\times 1)$ horizontal dominoes in $\lambda$ swap their placement (between the first two columns and the first row) without altering their orientation. We name the operation that describes this process the flip involution. It is a partial symmetry of $\SS^{(2)}[\SS^{(3)}[\SS^{(2)}[V]]]$ because it is only defined on the coefficients $a_{\ydiagram{2}[\ydiagram{3}[\ydiagram{2}]]}^{\lambda}$ for $\lambda$ a hook+column.
\end{ex}

This situation evokes  the partial symmetry for the Littlewood--Richardson coefficients, originally conjectured to exist by Pelletier and Ressayre \cite{pelletier}, and described explicitly and shown to hold by Grinberg \cite{grinberg}.
The presence of symmetries often gives us a better grasp of these coefficients. They can be useful in simplifying the number of cases that need to be addressed in both proofs and calculations involving them \cite{briand15}.\\

In Section \ref{sec:results}, we provide an overview of the necessary background on symmetric functions, and present the following results.
We introduce a partial involution, defined on hook+column partitions, in Definition \ref{def_flip}, named the flip involution.  Lemma \ref{flip_combinatorially}  shows that this involution can be understood as a ``transposition" of brick tilings of partitions. We then give closed expressions for the multiplicities of hook+column irreducible representations in  $\SS^{(c)}[\SS^{(b)}[\SS^{ (a)}[V]]]$, when either $b$ or $c$ is equal to 2, in Theorems \ref{sc o s2 o sa} and \ref{s2 o sb o sa}. As a corollary, these are examples of flip-symmetric representations.
In Theorem \ref{Symmetry},  we show that the Schur functors $\SS^{(2)}$ and $\SS^{(1,1)}$ preserve flip-symmetry. 
This result  allows us to produce  infinite families of flip-symmetric representations. In Section \ref{proof}, we furnish proofs for all these results.
In Section \ref{final_comments}, we finish this article presenting  an analysis of certain intriguing sequences constructed in a natural way from flip-symmetric iterated plethysm coefficients that appear to be both unimodal and asymptotically normal. 

\section{Background and presentation of the results}\label{sec:results}

Let $V$ a complex vector space of dimension $n.$
 In the language of symmetric functions, the role of the irreducible representation $\SS^{\lambda}[V]$ of $\GL(V)$ is played by
  the Schur polynomial $s_\lambda(x_1, x_2, \ldots, x_n).$
  A representation $W$ of $\GL(V)$  will then corresponds to a Schur positive
  symmetric function $f$, and the multiplicity of $\SS^{\lambda}[V]$
  in $W$ is equal to the coefficient of $s_\lambda$ in $f$.

We follow Stanley \cite{stanleyEC2} for the standard concepts and notations in the theory of symmetric functions, the main exception being that we write our partitions using the French convention \cite{macdonald}.
A \emph{partition} is a weakly decreasing sequence of natural numbers in which there are finitely many non-zero entries. The nonzero entries  are called the parts. Two partitions
are equal if they have the same parts. The weight of a partition is defined as the
sum of its parts.

\newcommand{\Par}{\textnormal{Par}}
We let $\Par(N)$ be the set of partitions of weight $N$. We identify partitions with their Young diagrams.  The Young diagram of a partition $\lambda = (\lambda_1, \lambda_2, \ldots)$ is the set $\{(c,r)\ :\ 0\le c\le \lambda_r\}$, whose elements we call cells. According to the
French convention Young diagrams are drawn bottom-left justified, thus cells of the diagram are described by its Cartesian coordinates.

The transpose $\lambda'$ of $\lambda$ is the partition whose diagram is the image of the diagram of $\lambda$ under the reflection $(c,r)\mapsto(r,c)$. Given partitions $\mu$ and $\lambda$, the skew-partition $\mu/\lambda$ is the set of cells in $\mu$ but not in $\lambda$. We define the sum of two partitions $\lambda$ and  $\mu$ as the partition $\lambda+\mu = (\lambda_1+\mu_1, \lambda_2+\mu_2,\ldots)$. The union of two partitions is the sorting of its parts.

\begin{de}
    We say a partition is a \emph{hook+column partition} if it is of the form $(\alpha, 2^\beta, 1^\gamma)$, that is, the sum of a hook partition and a column partition.
\end{de}
Hook+column partitions were introduced by Langley and Remmel \cite{langley} along with hook+row partitions $(\alpha,\beta,1^\gamma)$, which are the union of a hook and a row partition.

\begin{warning}
When we write $\lambda = (\lambda_1, 2^{m_2(\lambda)}, 1^{m_1(\lambda)})$ for a hook+column partition, we   use the  convention that $\lambda_1$ can be equal to either one or two.
We write $(1^5)$ as $(1, 2^0, 1^4)$, and  $(2^2,1)$ as $(2,2^1,1^1).$
Thus, $m_1(\lambda)$ and $m_2(\lambda)$ are not always the multiplicities of 1 and 2 in $\lambda$.
\end{warning}

The algebra of symmetric functions $\Lambda$ is defined as the algebra $\Q[p_1, p_2, \ldots]$ generated by  a set of variables that play the role of  the \emph{power sum symmetric functions}. Note that these are not defined in terms of another set of variables.  

However, it will sometimes be useful to identify an element $f$ of $\Lambda$ with a  formal power series.  For this, we let $X=x_1+x_2+\cdots$ be an alphabet. (The $x_i$'s are variables.)  Then, we identify  $f\in\Lambda$ with its image $f[X]$ under the algebra morphism that maps $p_k$ to $x_1^k + x_2^k +\cdots$.  We write $f[X] = f(x_1, x_2, \ldots)$ and say that it is the evaluation of $f$ in $X$. In particular, we identify $p_1$ with $X=x_1+x_2+\cdots$.

The notion of \emph{plethysm of symmetric functions} comes from that of composition of representations.  Consider the action of a Schur functor $\SS^\lambda$ on a diagonal endomorphism whose eigenvalues are  variables. Its trace  is then the Schur polynomial $s_\lambda$. Therefore to compute the plethysm of two symmetric functions one has to substitute the monomials appearing in the first symmetric function into the second  one.
More precisely, let $f$ and $g$ be elements of $\Lambda$. If $g[X]$ is a sum of monic monomials, 
$g[X] = g_1 + g_2 + \cdots$ then 
$
\big( f \circ g \big) [X]=f(g_1, g_2, \ldots)
$.

\begin{ex}\label{alphabet}
Let $f[X]= 2p_2[X]$, then, since $2p_2[X] = 2x_1^2 + 2x_2^2 +\cdots = x_1^2 + x_1^2 + x_2^2 + x_2^2 + \cdots$, we have that
\[
p_n\, [ \, 2p_2 [X] \,] = p_n(x_1^2, x_1^2, x_2^2, x_2^2,\ldots) = 2\, p_{2n}[X].
\]
\end{ex}

\begin{warning}\label{note: different evaluations}
One needs to be careful using plethystic notation. In general, evaluating on alphabet $cX$ is not equivalent to evaluating on the alphabet $cx_1 + cx_2 + \cdots$. We denote the first by $f[cX]$ and the latter by $f[tX]|_{t=c}$. In particular, $-p_k[X] = p_k[-X] \ne \left.p_k[tX]\right|_{t=-1} = (-1)^k p_k[X]$.
\end{warning}

The plethysm of symmetric functions can be defined axiomatically.

\begin{samepage}
\begin{de}[Plethysm of symmetric functions] The plethysm of symmetric functions is the operation $\circ : \Lambda\times\Lambda\to\Lambda$ verifying
\begin{enumerate}
    \item $p_n\circ p_m = p_{nm}$ for all $n, m\in\N_0$.
    \item For any $f\in\Lambda$, the map $g\mapsto g\circ f$ is a $\Z$-algebra homomorphism on $\Lambda$.
    \item For any $f\in\Lambda$, the equality $p_n\circ f = f\circ p_n$ holds.
\end{enumerate}
\end{de}
\end{samepage}

The core tools of this work come from plethystic calculus, namely, from the operation of evaluation in sums and differences of alphabets.  The following lemma is a consequence of the well known expansion of $s_\lambda\circ(f\pm g)$   found in \cite{macdonald}.
Let $c_{\mu,\nu}^\lambda$ denote the Littlewood--Richardson 
coefficient indexed by partitions $\mu, \nu$ and $\lambda.$ 
That is, $c_{\mu,\nu}^\lambda$ is the coefficient of $s_\lambda$ in
the product $s_\mu\cdot s_\nu.$

\begin{figure}[h]
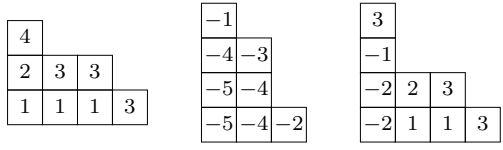

    \centering
    \ytableaubig
    \footnotesize
    \ytableaushort{4,233,1113}$\quad\quad$
    \ytableaushort{{-1},{-4}{-3},{-5}{-4},{-5}{-4}{-2}}$\quad\quad$
    \ytableaushort{3,{-1},{-2}23,{-2}113}
    \ytableausetup{nosmalltableaux}
    \caption{Three valid SSYT with positive and/or negative letters.}
    \label{fig:negative}
\end{figure}

\begin{lem}\label{sergeev} Let $X$ and $Y$ be two alphabets and let $\lambda$ be a partition. Then:
\begin{enumerate}
    \item $s_\lambda[-X] = (-1)^{|\lambda|}s_{\lambda'}[X].$
    \item $s_\lambda[X+Y] = \sum_{\mu\subset\lambda} \ s_\mu[X] \ \cdot \ s_{\lambda/\mu}[Y] = \sum_{\mu,\nu}c_{\mu,\nu}^\lambda \ s_\mu[X]\ \cdot \ s_{\nu}[Y].$
    \item
    \( s_\lambda[X-Y] = \sum_{\mu\subset\lambda}(-1)^{|\lambda/\mu|}\ s_\mu[X]\ \cdot \ s_{(\lambda/\mu)'}[Y]
    \)
    \\ \phantom{$s_\lambda[X-Y]$} 
    \(
    = \sum_{\mu,\nu}(-1)^{|\nu|}\ c_{\mu,\nu}^\lambda \ s_\mu[X]\ \cdot \ s_{\nu'}[Y].
    \)
\end{enumerate}
\end{lem}
\begin{Note}\label{note:positive and negative} Let $X$ and $Y$ be two alphabets. Then, Lemma \ref{sergeev} says that $s_\lambda[X-Y]$ is the generating function of the tableaux $T$ on positive letters from $X$ and negative letters from $-Y$ obeying the semistandard rules for the positive entries and a similar, but opposite rule for the negative ones: Negative entries should be weakly increasing across the columns when reading from from bottom to top, and strictly increasing across the rows, when reading from left to right. See Figure \ref{fig:negative}.
\end{Note}

\begin{de}\label{def_flip} 
Fix a non-negative number $\gamma$, and
let $f$ be a Schur positive and homogeneous symmetric function of
degree $n$.  Define $b_{f,\gamma}^\beta$ as the coefficient  of $s_{(n-2\beta-\gamma, 2^\beta, 1^\gamma)}$ in $f$. Then 
The \emph{hook+column sequence}  $\Sigma(f, \gamma)$ is defined as the sequence
\[
\Sigma(f, \gamma) = (b_{f,\gamma}^\beta)_{\beta \ge 0}.
\]
We refer to a sequence as \emph{symmetric} if its non-zero entries form a symmetric sequence.
\end{de}

\begin{ex} Tables \ref{tab:snos3os2} and \ref{tab:snos4os2} give examples of hook+column sequences. The data clearly indicates that these are symmetric sequences. Some further properties of these sequences will be discussed in the final section of this work.
\end{ex}

\begin{table}[h]
    \centering
    {
    \footnotesize 
    \begin{tabular}{  l | l | l }
    \toprule
    Function, $f$ & $\gamma$  & Hook+column sequence, $\Sigma(f, \gamma)$  \\
    \midrule
    $s_1\circ s_3\circ s_2$ & 0 & (1, 1, 1)\\
    $s_2\circ s_3\circ s_2$ & 0 & (1, 2, 3, 3, 2, 1)\\
    $s_3\circ s_3 \circ s_2$ & 0 & (1, 2, 5, 7, 8, 7, 5, 2, 1)\\
    $s_4\circ s_3 \circ s_2$ & 0 & (1, 2, 5, 10, 15, 18, 18, 15, 10, 5, 2, 1)\\
    $s_5\circ s_3 \circ s_2$ & 0 & (1, 2, 5, 10, 15, 19, 28, 36, 38, 36, 28, 19, 10, 5, 2, 1)\\
    $s_6\circ s_3 \circ s_2$ & 0 & (1, 2, 5, 10, 19,  33, 49, 63, 72, 72, 63, 49, 33, 19, 10, 5, 2, 1)\\
    \bottomrule
    \end{tabular}}
    \vspace{.5em } 
    \caption{The hook+column sequence of $s_c\circ s_b\circ s_2$ and $\gamma=0$, for $b= 3$ and $c=1, \ldots, 6$.}
    \label{tab:snos3os2}
\end{table}

\begin{table}[h]
    \centering
    {
    \footnotesize 
    \begin{tabular}{  l | l | l }
    \toprule
    Function, $f$ & $\gamma$  & Hook+column sequence, $\Sigma(f, \gamma)$  \\
    \midrule
    $s_1\circ s_4\circ s_2$ & 0 & (1, 1, 1, 1)\\
    $s_2\circ s_4\circ s_2$ & 0 & (1, 2, 3, 4, 4, 3, 2, 1)\\
    $s_3\circ s_4 \circ s_2$ & 0 & (1, 2, 5, 8, 11, 13, 13, 11, 8, 5, 2, 1)\\
    $s_4\circ s_4 \circ s_2$ & 0 & (1, 2, 5, 11, 18, 26, 34, 38, 38, 34, 26, 18, 11, 5, 2, 1)\\
    $s_5\circ s_4 \circ s_2$ & 0 & (1, 2, 5, 11, 22, 36, 55, 74, 90, 100, 100, 90, 74, 55, 36, 
    \ldots) 
    \\
    $s_6\circ s_4 \circ s_2$ & 0 & (1, 2, 5, 11, 22, 41, 68, 103, 144, 184, 217, 236, 236, 217, 184,
    \ldots) 
    \\
    \bottomrule
    \end{tabular}}
    \vspace{.5em } 
    \caption{The hook+column sequence of $s_c\circ s_b\circ s_2$ and $\gamma=0$, for $b= 4$ and $c=1, \ldots, 6$.}
    \label{tab:snos4os2}
\end{table}

We describe the symmetry  present in these sequences using the flip involution, a partial involution on partitions that we proceed to define.

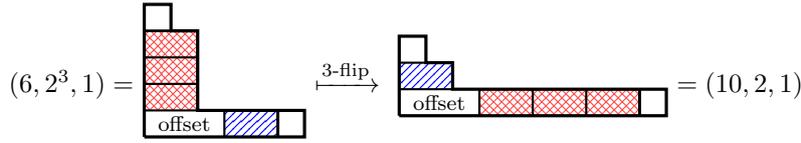
\begin{figure}[H]
    \centering
    $(6,2^3,1) =$
    \begin{tikzpicture}[x = 1em, y = 1em, baseline=1.7em]
        \fill[pattern = north east lines, pattern color = blue]
        (3,0) -- (3,1) -- (5,1) -- (5,0);
        \fill[pattern = crosshatch, pattern color = red!70]
        (0,1) -- (0,4) -- (2,4) -- (2,1);
        \draw[very thick]
        (0,0) -- (6,0) -- (6,1) -- (2,1) -- (2,4) -- (1,4) -- (1,5) -- (0,5) -- (0,0);
        \draw[thick]
        (3,0) -- (3,1) -- (5,1) -- (5,0)
        (0,1) -- (0,4) -- (2,4) -- (2,1)
        (0,1)--(2,1)
        (0,2)--(2,2)
        (0,3)--(2,3)
        (0,4)--(2,4)
        (5,0)--(5,1);
        \filldraw
        (1.5,0.5) node{\footnotesize offset};
    \end{tikzpicture}$\ \xmapsto{\text{3-flip}}\ $
    \begin{tikzpicture}[x = 1em, y = 1em, baseline=0.9em]
        \fill[pattern = crosshatch, pattern color = red!70]
        (3,0) -- (3,1) -- (9,1) -- (9,0);
        \fill[pattern = north east lines, pattern color = blue]
        (0,1) -- (2,1) -- (2,2) -- (0,2);
        \draw[very thick] 
        (0,0) -- (10,0) -- (10,1) -- (2,1) -- (2,2) -- (1,2) -- (1,3) -- (0,3) -- (0,0);
        \draw[thick]
        (3,0) -- (3,1) -- (9,1) -- (9,0)
        (0,1) -- (2,1) -- (2,2) -- (0,2)
        (3,0)--(3,1)
        (5,0)--(5,1)
        (7,0)--(7,1)
        (9,0)--(9,1);
        \filldraw
        (1.5,0.5) node{\footnotesize offset};
    \end{tikzpicture}
    $= (10,2,1)$
    \\
    \caption{The 3-flip of  $(6,2^3,1)$, is
    $\flip(3;(6,2^3,1))=(10, 2,1)$.}
    \label{fig:flip1}
\end{figure}

\begin{de} 
Let $\lambda$ be a hook+column partition, and let $r\ge2$ be an integer.
Assume furthermore that one can write $\lambda_1 - r - \gamma = 2\delta$ for some non-negative integer $\delta$.
The \emph{flip involution with offset $r$}, also called $r$-flip, is defined as follows:
\[
\lambda = (r + 2\delta + \gamma, 2^\beta, 1^\gamma) 
\mapsto
\flip(r; \lambda) = (r + 2\beta + \gamma, 2^\delta, 1^\gamma).
\]
Note that $\flip(r; -)$ is clearly an involution on its domain. See Figure \ref{fig:flip1}.
\end{de}

The following lemma translates the $r$-flip involution to a tiled transposition of Young diagrams. We say a partition is tiled if we have a collection of non-overlapping rectangles (tiles) of shape $(h\times w)$ for some $h,w\in\N$ covering its Young diagram.
A tiling is a brick tiling if every tile is of height 1.

A brick tiling of a Young diagram corresponds to a tableau in the following manner: A row tiled into $(1\times t_1), (1\times t_2),\ldots, (1\times t_k)$ can be collapsed to the row-tableau $\ytableaux \footnotesize\ytableaushort{{t_1}{t_2}{\cdots}{t_k}}\,$; now the diagram can be collapsed row by row into a tableau $T_\lambda$. (The shape of this tableau need not be a partition, but a composition.) If the shape of $T_\lambda$ is a partition, we can define the tiled-transpose of $\lambda$ as the (brick tiled) partition $\mu$ whose tableau $T_\mu$ is the transpose of $T_\lambda$.

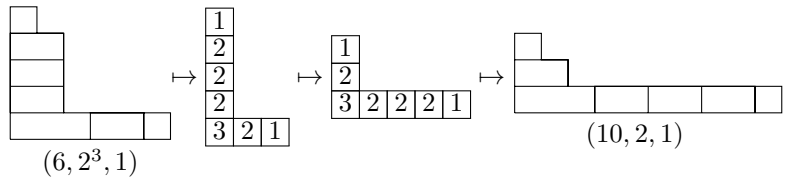
\begin{figure}[H]
    \centering
    \begin{tikzpicture}[x = 1em, y = 1em, baseline=2.1em]
        \draw
        (3,0) -- (3,1) -- (5,1) -- (5,0)
        (0,1) -- (0,4) -- (2,4) -- (2,1)
        (0,0) -- (6,0) -- (6,1) -- (2,1) -- (2,4) -- (1,4) -- (1,5) -- (0,5) -- (0,0)
        (0,1)--(2,1)
        (0,2)--(2,2)
        (0,3)--(2,3)
        (0,4)--(2,4)
        (5,0)--(5,1)
        (3,0) node[anchor=north] {$(6,2^3,1)$};
    \end{tikzpicture}$\mapsto
    \ytableaux \ytableaushort{1,2,2,2,321} \mapsto
    \ytableaux \ytableaushort{1,2,32221} \mapsto$
    \begin{tikzpicture}[x = 1em, y = 1em, baseline=1.1em]
        \draw
        (3,0) -- (3,1) -- (9,1) -- (9,0)
        (0,1) -- (2,1) -- (2,2) -- (0,2)
        (0,0) -- (10,0) -- (10,1) -- (2,1) -- (2,2) -- (1,2) -- (1,3) -- (0,3) -- (0,0)
        (3,0)--(3,1)
        (5,0)--(5,1)
        (7,0)--(7,1)
        (9,0)--(9,1)
        (4.5,0) node[anchor=north] {$(10,2,1)$};
    \end{tikzpicture}
    \\
    \caption{The tiled transpose of this tiling of $(6,2^3,1)$ is
    $\flip(3;(6,2^3,1))=(10, 2,1)$.}
    \label{fig:r-flip pos}
\end{figure}

\begin{lem}\label{flip_combinatorially}
Let $\lambda = (\lambda_1, 2^{\beta}, 1^{\gamma})$ be a hook+column partition and let $r\ge2$ be an integer. Assume furthermore that one can write $\lambda_1 - r - \gamma = 2\delta$ for some non-negative integer $\delta$.
Consider the following brick tiling of $\lambda$.
    \begin{enumerate}
        \item The first row is made up of one  $(1\times r)$-tile, followed by $\delta$ tiles of shape $(1\times 2)$, and $\gamma$ tiles of shape $(1\times 1)$.
        \item Each of the remaining rows forms a tile.
    \end{enumerate} 
Then the flip involution with offset $r$ corresponds to tiled transposition. See Figure \ref{fig:r-flip pos}.
\end{lem}

Given a symmetric function $f$ and a partition $\mu$, we denote  the coefficient of $s_\mu$ in the expansion of $f$ in the Schur basis by $[\mu]\,f$. We define the \emph{support} of $f$ as the set of partitions appearing in the decomposition of $f$ in the Schur basis, $\supp(f)=\{\mu \ : \ [\mu]
\, f\ne0\}.$

\begin{de}
    A symmetric function $f$ is \emph{flip-symmetric with offset $r$} if for all hook+column $\mu$ in $\supp(f)$,
    \begin{enumerate}
    \item $\flip(r;-)$ is defined on $\mu$, and \item $[\nu]f = [\flip(r;\mu)]f$.
    \end{enumerate} Moreover, $f$ is \emph{flip-symmetric} if its flip-symmetric with some offset $r$.
\end{de}

In particular, hook+column sequences of flip-symmetric functions are symmetric sequences.
However, as the following example illustrates, not all symmetric functions with symmetric hook+column sequences are flip-symmetric.

\begin{ex}\label{ex:counterex}
Here are some functions with symmetric hook+column sequences arising from plethysm.
\begin{align*}
  &s_{1,1}\circ s_{1,1} = s_{(2, 1, 1)}
  && \text{Not flip-symmetric.} \\
  &s_{1,1}\circ s_{2} = s_{(3, 1)}
  &&\text{Flip-symmetric with offset $2$.}\\
  &s_{2}\circ s_{1,1} = s_{(1^4)}+ s_{(2,2)}
  && \text{Not flip-symmetric.} \\
  &s_{2}\circ s_{2} = s_{(4)} + s_{(2,2)}
  &&\text{Flip-symmetric with offset $2$.}
\end{align*}
For instance, the third is not flip-symmetric because there is no $r$ such that $\flip(r;-)$ fixes the set $\{(1^4), (2^2)\}$.
\end{ex}

This observation calls for a description of all flip-symmetric iterated plethysms. The following formula of Langley and Remmel 
will provide us with the first non-trivial examples of such functions, and let us construct two more flip-symmetric families of functions arising from plethysm.
Given a symmetric function $f = \sum a_\lambda s_\lambda$, we introduce the notation $\hc{f}$ for $\sum_{\lambda = (\alpha,\beta,1^\gamma)} a_\lambda s_\lambda$. If furthermore we write $\hcG{f}$, we are fixing $\gamma$  for all $\lambda$.
\begin{thm}[Langley and Remmel \cite{langley}, Thm. 4.8]\label{sasb} For any $a,b\ge2$, 
\[\hc{s_b\circ s_a} = \hcZ{s_b\circ s_a} = \sum_{k<b}s_{(ab-2k,2^k)}.
\]
\end{thm}

We generalize this theorem of Langley and Remmel, and obtain closed formulas for the iterated plethystic coefficients $a_{(c)[(b)[(a)]]]}^{\lambda}$ when either $b$ or $c$ is equal to 2, and $\lambda$ is a hook+column partition. As a 
corollary, we obtain that both families of iterated plethysms are flip-symmetric.

\begin{thm}\label{s2 o sb o sa}
   Let $a$ and $b$ be integers $\ge2$. Then,
    \begin{align*}
    \hc{s_2\circ s_b\circ s_a} =  &\sum\limits_{k=0}^{2b-1}  \min\left\{k+1, 2b-k\right\} \cdot s_{(2ab-2k,\ 2^{k})}\\
    &+\sum\limits_{k=1}^{2b-3} \min\left\{\left\lfloor\frac{k+1}{2}\right\rfloor, \left\lfloor\frac{2b-1-k}{2}\right\rfloor\right\} \cdot s_{(2ab-2k-1,\ 2^{k},\ 1)}.
     \end{align*}
\end{thm}

The proof of this theorem can be found in Section \ref{sec:s2sasb} and relays on Langley and Remmel's formula. It reduces our calculation to counting the number of integer points within particular polytopes.

   \begin{thm}   \label{sc o s2 o sa}
   Let $a$ and $c$ be integers $\ge2$. Then, $    \hc{s_c\circ s_2\circ s_a} $ is equal to
   \begin{align*}
 &\sum\limits_{k=0}^{2c-1}\min\left\{\frac{k^2+k+2}{2},\frac{(2c-1-k)^2+(2c-1-k)+2}{2}\right\}\cdot s_{(2ac-2k, 2^k)}\\
    &\qquad\qquad+ \sum\limits_{k=1}^{2c-3}\min\left\{\left\lfloor\frac{(k+1)^2}{4}\right\rfloor,\left\lfloor\frac{(2c-1-k)^2}{4}\right\rfloor\right\}\cdot s_{(2ac-2k-1,\ 2^k,\ 1)}.
   \end{align*}

\end{thm}

The proof of this theorem appears in Section \ref{scs2sa}. The techniques used in this proof originate in \cite{rosas}, and \cite{langley}. They relay on careful evaluation of $s_c\circ s_2\circ s_a$ on signed finite alphabets.

\begin{cor}
The iterated plethysms $s_c\circ s_b\circ s_a$ are flip-symmetric with offset $abc-2(bc-1)$ when either $b$ or $c$ equals 2.
\end{cor}

Note that the increasing parts of of the sequences of coefficients
appearing in Theorems \ref{s2 o sb o sa} and  \ref{sc o s2 o sa} 
are, respectively,  the  natural numbers, the natural numbers (repeated
twice), the triangular numbers plus one, also known as the
central polygonal numbers (OEIS A000124), and the ``quarter-squares'' (OEIS A002620 \cite{oeis}). The generating functions for all these sequences are rational.

\begin{Note}  There exists a relationship between hook+column sequences and the analogously defined hook+row sequences, via the $\omega$ involution \cite{macdonald}: for two homogeneous symmetric functions, $\omega(f\circ g) = f\circ \omega(g)$ if degree($g$) is even and $\omega(f\circ g) = \omega(f)\circ \omega(g)$ otherwise. Consequently, if $\hc{s_{n_1}\circ s_{n_2}\circ\cdots\circ s_{n_k}}^\gamma$ $=\sum_\beta a_\beta s_{\nu_\beta}$, then applying the aforementioned formulas successively will give a function on the left-hand-side whose hook+row sequences coincide with the hook+column sequences of the original function. Our work could be therefore translated to hook+row sequences.
\end{Note}

We show that
the plethysm operation with either $s_{(2)}$ or $s_{(1,1)}$ preserves the flip-symmetry of Schur positive symmetric functions, a result that allows us to construct even more families of flip-symmetric functions.

\begin{thm}\label{Symmetry} Let $f\in\Lambda_n$ be a Schur positive symmetric function.
If $f$ is flip-symmetric with offset $r$, then both $s_{1,1}\circ f$ and $s_2\circ f$ are Schur positive flip-symmetric (with offset $2r-2$) symmetric functions in $\Lambda_{2n}$.
\end{thm}

On Table \ref{fig:normalidad}, we tabulate the hook+column sequences obtained from the iterated plethysms $s_2^{\circ k}$ for $\gamma=0$ and $k=2,3,4,5$.

As a corollary of Theorem \ref{Symmetry}, we obtain that the plethystic action of $h_{1,1} = s_2 + s_{1,1}$ also preserves flip-symmetry and Schur positivity.
In the process, we show that $(p_{1,1}\circ f)$ and $(p_2\circ f)$ also preserve flip-symmetry (although $p_2$ does not, in general, preserve Schur positivity). 

\begin{cor}
Let $f$ be a Schur positive homogeneous symmetric function $f$, and let $\lambda^1$, $\lambda^2$, $\ldots$,
$\lambda^k$
be any sequence of partitions of two. If $f$ is flip-symmetric with offset $r$, then  
\[s_{\lambda^1}\circ s_{\lambda^2}\circ\ldots\circ s_{\lambda^k}\circ f,\]
is also flip-symmetric with offset $2^kr - 2^{k+1}+2$.
\end{cor}
See Figures \ref{fig:arrows} and \ref{fig:normalidad} for some examples of flip-symmetric 
sequences obtained in this way.

\begin{figure}[H]
\begin{equation*}
    \ytext
    \hc{s_{2}\circ s_{2} \circ s_{3}} =
    \tikzmark{6222}{
        s_{\begin{tikzpicture}[x=0.35em, y=0.35em, baseline=1.2em]
            \draw
            (0,1)--(0,2)--(2,2)--(2,1)
            (0,2)--(0,3)--(2,3)--(2,2)
            (0,3)--(0,4)--(2,4)--(2,3)
            (0,0)--(6,0)--(6,1)--(0,1)--(0,0);
        \end{tikzpicture}}
    } + 
    \tikzmark{822}{
        2s_{\begin{tikzpicture}[x=0.35em, y=0.35em, baseline=0.85em]
            \draw
            (0,1)--(0,2)--(2,2)--(2,1)
            (0,2)--(0,3)--(2,3)--(2,2)
            (6,0)--(8,0)--(8,1)--(6,1)
            (0,0)--(6,0)--(6,1)--(0,1)--(0,0);
        \end{tikzpicture}}
    } +
    \tikzmark{102}{
        2s_{\begin{tikzpicture}[x=0.35em, y=0.35em, baseline=0.5em]
            \draw
            (0,1)--(0,2)--(2,2)--(2,1)
            (8,0)--(10,0)--(10,1)--(8,1)
            (6,0)--(8,0)--(8,1)--(6,1)
            (0,0)--(6,0)--(6,1)--(0,1)--(0,0);
        \end{tikzpicture}}
    } + 
    \tikzmark{12}{
        s_{\begin{tikzpicture}[x=0.35em, y=0.35em, baseline=0.15em]
            \draw
            (10,0)--(12,0)--(12,1)--(10,1)
            (8,0)--(10,0)--(10,1)--(8,1)
            (6,0)--(8,0)--(8,1)--(6,1)
            (0,0)--(6,0)--(6,1)--(0,1)--(0,0);
        \end{tikzpicture}}
    } + 
    \tikzmark{921}{
        s_{\begin{tikzpicture}[x=0.35em, y=0.35em, baseline=0.85em]
            \draw
            (0,2)--(0,3)--(1,3)--(1,2)
            (0,1)--(0,2)--(2,2)--(2,1)
            (8,0)--(9,0)--(9,1)--(8,1)
            (6,0)--(8,0)--(8,1)--(6,1)
            (0,0)--(6,0)--(6,1)--(0,1)--(0,0);
        \end{tikzpicture}}
    }.
\end{equation*}
    \centering
    \begin{tikzpicture}[overlay, remember picture]
        \draw[\tip-\tip, thick,rounded corners] (6222.south)--($(6222) + (0,-1) $)--($(12) + (0,-1.19)$)--(12.south);
        \draw[\tip-\tip, thick,rounded corners] (822.south)--($(822) + (0,-0.83)$)--($(102)+(0,-0.9)$)--(102.south);
        \draw[-\tip, thick] ($(921) + (0,-0.5)$) arc (110:405:.25) -- ++(135:5pt);
    \end{tikzpicture}
    \vspace{.5em} 
    \caption{The function $s_2\circ s_2\circ s_3$ is flip-symmetric with offset 6. The arrows indicate the images of the 6-flip involution $\flip(6; -)$. The indexing partitions appear tiled according to the flip algorithm.}
    \label{fig:arrows}
\end{figure}

\section{Proofs} \label{proof}

\subsection{A handful of lemmas}\label{sec:handful}
\ytableausetup{boxsize=0.5em}

Starting from the elementary remark that $2 s_2 = p_{1,1} + p_2$ and $2 s_{1,1} = p_{1,1} - p_2$, we get that the plethysms $s_2\circ f$ and $s_{1,1}\circ f$ are completely determined by the plethysms $p_{2}\circ f$ and $p_{1,1}\circ f$. The following lemma says that when we want to compute $\hc{s_\sigma\circ f} $, we can restrict our attention 
to  plethysms of the form $p_{2}\circ\hc{f}$ and $p_{1,1}\circ\hc{f}$. 

\begin{lem}\label{s2 o f} Let $f$ be a symmetric function. Then, for any partition $\sigma$, we have 
\[
\hc{s_\sigma\circ f} = \hc{s_\sigma \circ \hc{f}}.
\]
\end{lem}
 
\begin{proof}
Let $f = \sum_\mu d_\mu s_\mu$ be a symmetric function. We follow the outline of the proof of \cite[Theorem 5.1]{gutierrez}. For simplicity, we also use the same notation. We have
\begin{align*}
    s_\sigma\circ f
    &= \sum_\lambda \frac{\chi^\sigma(\lambda)}{z_\lambda} p_\lambda\circ f
    = \sum_\lambda \frac{\chi^\sigma(\lambda)}{z_\lambda} \prod_i p_{\lambda_i}\circ f
    = \sum_\lambda \frac{\chi^\sigma(\lambda)}{z_\lambda} \prod_i f\circ p_{\lambda_i}\\
    &= \sum_\lambda \frac{\chi^\sigma(\lambda)}{z_\lambda} \prod_i \sum_\mu d_\mu s_\mu \circ p_{\lambda_i}
    = \sum_\lambda \frac{\chi^\sigma(\lambda)}{z_\lambda} \prod_i \sum_\mu d_\mu \sum_\tau b_{\lambda_i,\mu}^\tau s_\tau.
\end{align*}
The expansion $s_\sigma= \sum_\lambda \frac{\chi^\sigma(\lambda)}{z_\lambda} p_\lambda$ is classical. The coefficients $b_{\lambda_i,\mu}^\tau$ appearing the the last equality come from the expansion $s_\mu \circ p_{\lambda_i} = \sum_\tau b_{\lambda_i,\mu}^\tau s_\tau$.

By the Littlewood--Richardson rule, a hook+column partition is a constituent of a product $s_\tau\cdot s_\pi$ only if $\tau$ and $\pi$ are hook+column partitions.
From \cite[Theorem 4.1]{gutierrez}, we know $b_{\lambda_i,\mu}^\tau$ is nonzero only if $\mu\subseteq\tau$ (that is, $\mu_i\le\tau_i$ for all $i$). Then,
\begin{align*}
    \hc{s_\sigma\circ f}
    &= \hc{\sum_\lambda\frac{\chi^\sigma(\lambda)}{z_\lambda} \prod_i \sum_\mu d_\mu \sum_{\tau} b_{\lambda_i,\mu}^\tau s_\tau}\\
    &= \hc{\sum_\lambda\frac{\chi^\sigma(\lambda)}{z_\lambda} \prod_i \sum_\mu d_\mu \sum_{\tau \text{ h+c}} b_{\lambda_i,\mu}^\tau s_\tau}\\
    &= \hc{\sum_\lambda\frac{\chi^\sigma(\lambda)}{z_\lambda} \prod_i \sum_{\mu \text{ h+c}} d_\mu \sum_{\tau \text{ h+c}} b_{\lambda_i,\mu}^\tau s_\tau} \\
    &= \hc{s_\sigma\circ\hc{f}}.
\end{align*}
\end{proof}

In \cite{carre}  Carré and Leclerc found an elegant description of the  plethystic action of $p_2$ in terms of domino tableaux and their $2$-signs. 
Let $\lambda$ be a partition of $2n$. A tiling is domino if every tile is $(1\times2)$ or $(2\times1)$.
A domino tableau of $\lambda$ is a labelling of a domino tiling of $\lambda$ with non-negative integers so that the numbers increase weakly in rows and strictly increase in columns.
There is a general definition for the $n$-sign of a partition (see \cite{wildon}). However, we are only interested in the $2$-sign of a partition, defined as
\[
    \sgn_2(\lambda):=(-1)^{\#\text{vertical dominoes in a domino tiling of $\lambda$}}.
\]
As a fact, this is independent of the domino tiling of the partition.

\begin{ex}
To compute the 2-sign of $(4,3,1)$, we first compute a domino tiling of its diagram.
\[
\ytableaux \lambda = \ydiagram{1, 3, 4} 
\hspace{2cm}
\begin{tikzpicture}[x=1em, y=1em, baseline=1.2em]
    \draw
        (0,0)--(4,0)
        (0,1)--(4,1)
        (1,2)--(3,2)
        (0,3)--(1,3)
        (0,0)--(0,3)
        (1,1)--(1,3)
        (2,0)--(2,1)
        (3,1)--(3,2)
        (4,0)--(4,1);
\end{tikzpicture}
\hspace{2cm}
\sgn_2(\lambda) = -1.
\]    
\end{ex}

Simple inspection leads to the following realization: if $\lambda=(\alpha,2^\beta,1^\gamma)$ is a hook+column partition, then the sign of $\lambda$ only depends on the congruence class of $\gamma$ modulo 4, as shown in Table \ref{tab:sign}.
\begin{table}[h]
    \centering\scriptsize
    \begin{tabular}{  l | l l l l }
    \toprule
$\gamma$ modulo 4 &0&1&2&3\\
\midrule
$\sgn_2(\mu)$&+&--&--&+\\
\bottomrule
\end{tabular}
\vspace{.5em } 
    \caption{The 2-sign of $\lambda=(\alpha,2^\beta,1^\gamma)$ only depends on $\gamma$.}
    \label{tab:sign}
\end{table}

Carré and Leclerc's description gives the following formula.

\begin{lem}\label{p2slambda} Let $\lambda = (\lambda_1, 2^{m_2(\lambda)}, 1^{m_1(\lambda)})$ be a hook+column. Then, $(-1)^{m_1(\lambda)} \hc{p_2\circ s_\lambda} $ is equal to
\begin{align*}
\sum_{\small{\beta\ge2m_2(\lambda)}} (-1)^{\beta} s_{(2\lambda_1,2^\beta,1^\bullet)}
- s_{(2\lambda_1-1, 2^{2m_2(\lambda)},1^\bullet)}
+ \sum_{\small{\beta\ge 2m_2(\lambda)+1}} (-1)^{\beta+1} s_{(2\lambda_1-2,2^\beta,1^\bullet)}.
\end{align*}
\end{lem}

There is an equivalent way of stating Lemma
\ref{p2slambda}. Given a symmetric function $f$ and a partition $\lambda$, let us denote by $[\lambda]f$ the coefficient of $s_\lambda$ in the expansion of $f$ in the Schur basis.
Let $\lambda = (\lambda_1, 2^\beta, 1^\gamma)$ be a hook+column partition, and let $f$ be a homogeneous symmetric function of degree $n$. Then,   
\[
[\lambda]\, (p_2\circ f) = \sgn_2(\lambda)\cdot\#\D{2}{},
\]
where $\D{2}{}\subseteq\mathtt{Par}(2n)$ is a multiset constructed as follows: 
\begin{enumerate}
    \item 
The underlying set is
    \[
    \begin{cases}
        \left\{\left(\frac{\lambda_1+1}{2},\ 2^{\frac{\beta}{2}},\ 1^{\frac{\gamma-1}{2}}\right)\right\}
        & \text{if } \lambda_1 \text{ is odd,}
        \\
        \phantom{.}
        \vspace{-0.75em} 
        \\
        \left\{\mu\ : \ \substack{2\mu_1=\lambda_1\\ m_2(\mu)\le\beta /2}\right\}\cup\left\{\mu\ : \ \substack{2\mu_1=\lambda_1+2\\ m_2(\mu)\le \frac{\beta-1}{2}}\right\}
        & \text{if } \lambda_1 \text{ is even.}
    \end{cases}
    \]
    
\item
    The multiplicity of  $\mu$ in $\D{2}{}$ is $[\mu] \, f$.
\end{enumerate}
\noindent
Observe that in the case when $\lambda_1$ is even, $\D{2}{}$ naturally splits into two multisets, which we call $\D{2}{0}$ and $\D{2}{2}$ respectively.\\

We need to  evaluate  the expression $\hc{p_
{1,1}\circ f}$. Since $p_{1,1}\circ f = f^2$, we just need  to compute $\hc{s_\lambda\cdot s_\mu}$ when both $\lambda$ and $\mu$ are hook+column partitions. The following formula comes from the Littlewood--Richardson rule.

\begin{lem}\label{p11slambda}
Let $\mu = (\mu_1, 2^{m_2(\mu)}, 1^{m_1(\mu)})$ and $\nu = (\nu_1, 2^{m_2(\nu)}, 1^{m_1(\nu)})$ be hook+columns.  Set $\alpha = \nu_1 + \mu_1$,  $m_2 =m_2(\mu) + m_2(\nu)$, and $m_1 = \min\{m_1(\mu), m_1(\nu)\}$. 

Then,
\begin{align*} 
\hc{s_\mu\cdot s_\nu} = &\sum\limits_{\beta = m_2}^{m_2 + m_1} s_{(\alpha, 2^\beta, 1^\bullet)}  
+ \sum\limits_{\beta=m_2}^{m_2 + m_1 + 1} \chi_{(\mu,\nu)}^\beta\cdot s_{(\alpha-1, 2^{\beta}, 1^\bullet)} \\
&+ \sum\limits_{\beta=m_2 +1}^{m_2 + m_1+1} s_{(\alpha-2, 2^\beta, 1^\bullet)},
\end{align*}
where $\chi_{(\mu,\nu)}^\beta=1$ if $\beta=m_2$ or $\beta=m_2+m_1+1$ and $\chi_{(\mu,\nu)}^\beta=2$ otherwise.
        \end{lem}

Consequently, if $\lambda = (\lambda_1, 2^\beta, 1^\gamma)$ and  $f$ is a homogeneous symmetric function of degree $n$, we can write $[\lambda] \, (p_{1,1}\circ f) = \#\D{1,1}{}$, where $\D{1,1}{}\subseteq\mathtt{Par}(2n)\times\mathtt{Par}(2n)$ is a multiset constructed as follows: the underlying set is the union of the sets
\begin{samepage}
    \begin{align*} 
    &\left\{(\mu,\nu)\ : \ \substack{\mu_1+\nu_1=\lambda_1\\ m_2\le\beta\le m_2+m_1}\right\},
    \qquad
    \left\{(\mu,\nu)\ : \ \substack{\mu_1+\nu_1=\lambda_1+1\\ m_2\le\beta\le 1+m_2+m_1}\right\}
    \\
    &\qquad \qquad \text{\qquad \qquad \qquad and} \qquad 
    \left\{(\mu,\nu)\ : \ \substack{\mu_1+\nu_1=\lambda_1+2\\ 1+m_2\le\beta\le 1+m_2+m_1}\right\}.
    \end{align*}
\end{samepage}
    (Therefore the multiset $\D{1,1}{}$ naturally splits into three multisets, which we call $\D{1,1}{0}$, $\D{1,1}{1}$ and $\D{1,1}{2}$ respectively.)
    The multiplicity of each pair $(\mu,\nu)$ in $\D{1,1}{0}$ and $\D{1,1}{2}$ is $[\mu]\, f\cdot [\nu] \, f$, whereas the multiplicity of $(\mu,\nu)$ in $\D{1,1}{1}$ is $\chi_{(\mu,\nu)}^\beta[\mu] \,f\cdot [\nu]\, f$.

\begin{Note}These multisets, or rather their underlying sets, can be interpreted as the sets of integral points of some polytopes in $\Z^6$, by identifying a hook+column pair $(\mu,\nu)=\big((\mu_1,2^{m_2(\mu)},1^{m_1(\mu)}),(\nu_1,$ $2^{m_2(\nu)},1^{m_1(\nu)})\big)$ with the point $\big(\mu_1,m_2(\mu),m_1(\mu),\nu_1,m_2(\nu),m_1(\nu)\big)$. Now, the equalities and inequalities that define our sets are viewed as the restriction to certain hyperplanes and regions of the space. Once this work is done, the integer points inside the intersection of those hyperplanes and regions form the announced polytope.
\end{Note}

Inspection of Lemmas \ref{p2slambda} and \ref{p11slambda} reveals a beautiful result. If $\mu = \nu$, then
\[\supp\!\left(\hc{p_2\circ s_\mu}\right)\subseteq\supp\!\left(\hcS{s_\mu^2}\right).\]
Furthermore, if $m_1(\mu)=0$, then the two sets are identical, and so are the multiplicities associated with each partition.
Finally, there is one more thing to notice: $\sgn_2((\alpha,2^\beta)) = 1$ for all $(\alpha, 2^\beta)$. Consequently, we obtain the following result.\begin{lem}\label{square} Let $\mu=(\alpha,2^\beta)$. Then, $\hcZ{p_2\circ s_\mu} = \hcZ{p_{1,1} \circ s_\mu}$.
\end{lem}

\subsection{An explicit formula for \texorpdfstring{$s_2\circ s_b\circ s_a$}{s2 o sb o sa} on hook+columns}\label{sec:s2sasb}

In this section, we prove Theorem \ref{s2 o sb o sa}. 
We build our proof on Langley and Remmel's Theorem  \ref{sasb}.
Note that their formula barely depends on $a$. In fact, all that $a$ introduces is a tail of zeros in our hook+column sequence. For instance, $\Sigma(s_3 \circ s_2, 0) = (1, 1, 1)$ and $\Sigma(s_3 \circ s_4,0) = (1, 1, 1, 0, 0, 0)$.
\begin{lem}\label{tail0}
The hook+column sequence $\Sigma(s_b\circ s_a,0)$ is $(1,\stackrel{b\text{ times}}{\ldots},1,0,\stackrel{Z-b\text{ times}}{\ldots},0)$, where $Z = \left\lfloor\frac{ab}{2}\right\rfloor$.
\end{lem}
\begin{proof}
    The hook+column partitions of size $ab$ with $\gamma=0$ are $(ab), (ab-2, 2), (ab-4, 2^2), \ldots, (ab-2(Z-1), 2^{Z-1})$. The claim now is an immediate consequence of Theorem \ref{sasb}.
\end{proof}
Consequently, we can suppose $a=2$ hereafter without loss of generality.

\begin{proof}[Proof of Thm.~\ref{s2 o sb o sa}]
We begin the proof observing that, since $s_2 = \frac{1}{2}(p_{2}+ p_{1,1})$, then, by Theorem \ref{sasb},
\begin{align}\label{eq:suma gorda}
    2s_2 \circ \hc{s_b\circ s_2}
    &= (p_2 + p_{1,1})\circ \sum_{k<b} s_{(2b-2k, 2^k)}\nonumber\\
    &= \sum_{k<b} p_2\circ s_{(2b-2k, 2^k)} + \sum_{k<b} s_{(2b-2k, 2^k)}^2 \nonumber\\
    &\qquad\qquad\qquad 
    + 2 \sum_{i<j<b} s_{(2b-2i, 2^i)}\cdot s_{(2b-2j, 2^j)}\nonumber\\
    &= \sum_{k<b} 2s_2\circ s_{(2b-2k, 2^k)} + 2 \sum_{i<j<b} s_{(2b-2i, 2^i)}\cdot s_{(2b-2j, 2^j)}.
\end{align}

We now restrict to the $\gamma=0$ part of the hook+column sequence.
Note that \cite[Thm.~4.8 (3)]{langley} gives $\hc{s_2\circ s_\lambda} = \hcZ{s_2\circ s_\lambda}$ whenever $\lambda = (\lambda_1, 2^\beta)$ is a hook+column with $\gamma=0$. With this and using Lemma \ref{square}, we simplify our expression into
\[
\hcZ{s_2\circ s_b \circ s_2} = \hcZ{\sum_{i\le j<b} s_{(2b-2i, 2^i)}\cdot s_{(2b-2j, 2^j)}}.
\]
We use Lemma \ref{p11slambda} to compute the products appearing in this equation. For each term in the second sum, we get the exactly the two hook+column partitions $\lambda$ with $m_1(\lambda)=0$,  whose first row verify
\begin{math}
\lambda_1 \in \left\{4b -2(i+j),\ 4b - 2(i+j+1)\right\}.
\end{math}
Consequently, the hook+column partitions that will appear in the sum are those whose first row is in the set $\left\{2, 4, 6, \ldots, 4b\right\}$. We now ask how many times  each one appear.

Take $b\in\N_{\ge2}$ and $k\in\left\{1,2,\ldots,2b\right\}$. We ask how many integer pairs $(i,j)$  there are in the polytope $\Delta := \left\{0\le i\le j<b\right\}$, which are solutions to \textit{either} of these two equations:
\begin{equation}\label{eq:sistema}
    \begin{cases}
4b-2(i+j) = 2k, \\
4b-2(i+j+1) = 2k,
\end{cases} \textnormal{ or, consequently, } \, \, \,
\begin{cases}
i+j=2b-k,\\i+j=2b-k-1.
\end{cases}
\end{equation}

Refer to Figure \ref{fig:politopos} for a graphical representation.

\begin{figure}[h]
    \centering
    \begin{tabular}{ccc}
    \begin{tikzpicture}[x=1em, y=1em]
    \ejes{5}{5}
    
    \filldraw[blue, opacity=0.1]
        (0,0)--(4,0)--(4,4);
    \draw[blue, thick]
        (-1,-1)--(6,6)
        (4,-1)--(4,6);
    
    \draw[red, thick]
        (-1,5)--(5,-1)
        (-1,6)--(6,-1);
        
    \draw[arrows=-{\tip},red,thick]
        (3,2)--(2.5,1.5);
    \draw[arrows=-{\tip},red,thick]
        (4,1)--(3.5,0.5);
        
    \fill[black]
        (2,2) circle (2pt)
        (3,2) circle (2pt)
        (3,1) circle (2pt)
        (4,1) circle (2pt)
        (4,0) circle (2pt);
        
    \fill[black]
        (6,0) node[anchor=west]{$j$}
        (0,6) node[anchor=south]{$i$}
        (6,6) node[anchor=north west]{$i=j$}
        (4,6) node[anchor=south]{$j=b-1$}
        (-1,5) node[anchor=east]{$i+j=2b-k-1$}
        (-1,6) node[anchor=east]{$i+j=2b-k$};
    \end{tikzpicture} &  $\qquad$ &
    \begin{tikzpicture}[x=1em, y=1em]
    \ejes{5}{5}
    
    \foreach \y in {0,...,3}
    \foreach \x in {\y,...,3}
    {
    \fill[black!60] (\x,\y) circle (1.2pt);
    \fill[black!60] (\x+0.5,\y+0.5) circle (1.2pt);
    }
    \foreach \y in {0,...,4}
    {
    \fill[black!60] (4,\y) circle (1.2pt);
    }
    
    \filldraw[blue, opacity=0.1]
        (0,0)--(4,0)--(4,4);
    \draw[blue, thick]
        (-1,-1)--(6,6)
        (4,-1)--(4,6);
    
    \draw[red, thick]
        (-1,5)--(5,-1);
    
    \fill[black]
        (2,2) circle (2pt)
        (2.5,1.5) circle (2pt)
        (3,1) circle (2pt)
        (3.5,0.5) circle (2pt)
        (4,0) circle (2pt);
    \end{tikzpicture}
    \end{tabular}
    \caption{We let $b=5$ and $k=5$. The polytope $\Delta$ is shaded in blue. Each black dot represents a valid pair. On the right, we illustrate the result of the described projection.}
    \label{fig:politopos}
\end{figure}
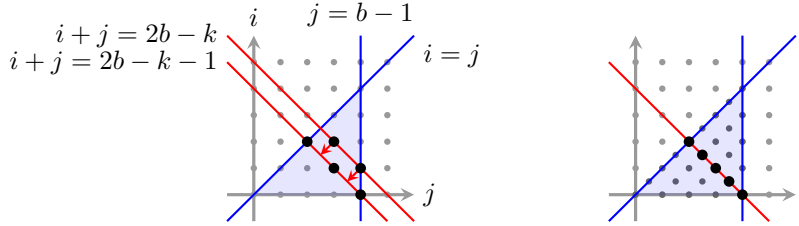

To help count the solutions, we will project them orthogonally from one of the lines to the other one, in such a way that all the projections remain inside $\Delta$. More precisely, if $k\ge b$ then project onto $\left\{i+j=2b-k\right\}$, an if $k<b$ to the other line.

One can easily see now what the coefficients are going to look like. Noting that the biggest line is counted twice (because we change the projection mid way), results in the desired integer sequence for $\gamma=0$.\\

Let us now bring back Equation \eqref{eq:suma gorda} and restrict to the $\gamma\ne0$ part of the hook+column sequence. Using \cite[Thm.~4.8 (3)]{langley} again, we obtain
\[
\hcN{s_2\circ s_b\circ s_2} = \hcN{\sum_{i<j<b} s_{(2b-2i,2^i)}\cdot s_{(2b-2j,2^j)}}.
\]
These products give rise, using Lemma \ref{p11slambda}, to a new polytope $\{0\le i< j< b\}$ together with the equation $i+j=2b-k$. A similar argument gives the desired sequence. This completes the proof of Theorem  \ref{s2 o sb o sa}.
\end{proof}

\subsection{An explicit formula for \texorpdfstring{$s_c\circ s_2\circ s_a$}{sc o s2 o sa} on hook+columns} \label{scs2sa}

This section is dedicated to the proof of Theorem \ref{sc o s2 o sa}. Fixing $\gamma=0$, we will show that the first $c$ terms of the hook+column sequence are given by the formula $\frac{n^2+n+2}{2}$ (OEIS A000124 \cite{oeis}) that have, as  generating function $\frac{z}{(1-z)^3}+\frac{1}{1-z}$ \cite{plouffe}. As before, it will be enough to show it for $a=2$. The techniques used in this section have been used by the second author in \cite{rosas} to study the Kronecker coefficients, and by Langley and Remmel \cite{langley}.

We begin by explicitly computing the evaluation of a Schur function indexed by a hook+column on the alphabet $1-x-y$.

\begin{lem}\label{slambda en 1-x-y} Let $\lambda=(\alpha, 2^\beta, 1^\gamma)$ and $\alpha\ge2$. Then,
\[
s_{\lambda}[1-x-y] =(-1)^\gamma (xy)^\beta (1-x)(1-y)\frac{x^{\gamma+1}-y^{\gamma+1}}{x-y}.
\]
In particular, the formula is independent of $\alpha$.
\end{lem}

\begin{proof} Recall Note \ref{note:positive and negative}. We construct all tableaux of shape $\lambda$ with three letters (fixing an order, let $1$ be the letter for the variable $1$, $-1$ for $-x$, and $-2$ for $-y$) and record the weight of each resulting tableau.

Notice that we have very little freedom when filling a hook+column with only these three letters. Our only choices are in the last entries in the first two columns (see Figure \ref{fig:slambda en 1-x-y}).
\begin{figure}[H]
    \centering
    \ytableaubig
    \footnotesize
    \ytableaushort{?,\vdots,??,{-2}{-1},\vdots\vdots,{-2}{-1},{-2}{-1}1\cdots1}
    \caption{Our only choices when filling a hook+column with 1, $-1$ and $-2$ are represented by a question mark.}
    \label{fig:slambda en 1-x-y}
\end{figure}

We can have 0, 1 or 2 entries equal to 1 in these cells. The rest can be filled with various quantities of $-1$s and $-2$s, resulting into weights
\[(xy - x - y + 1)\sum_{i+j=\gamma}(-1)^\gamma x^iy^j = (-1)^\gamma (1-x)(1-y)\frac{x^{\gamma+1}-y^{\gamma+1}}{x-y}.\]
As there are $\beta$ instances of \ytableaubig${
    \footnotesize\ytableaushort{{-2}{-1}}}$, the desired expression arises.
\end{proof}

Further inspection of Lemma \ref{slambda en 1-x-y} reveals that the restriction to $\gamma=0$ equates to restriction to monomials which are also monomials in the variable $(xy)$. We write this as follows:
\begin{equation}\label{eq:1}
\frac{\hcZ{s_c\circ s_2\circ s_2} [1-x-y]}{(1-x)(1-y)} = \left.\frac{s_c\circ s_2\circ s_2[1-x-y]}{(1-x)(1-y)}\right|_{(xy)}.
\end{equation}
The numerator $s_c\circ s_2\circ s_2[1-x-y]$ can be rewritten, using the equality $s_c = \sum_{\lambda\vdash c} z_\lambda^{-1} p_\lambda$ and the properties of plethysm, as
\begin{align*}
    s_c\circ s_2\circ s_2[1-x-y] 
    & = s_c \circ (s_4 + s_{2,2})[1-x-y] 
    && \text{(by \ref{sasb})}\\
    & = s_c [(1+xy)(1-x)(1-y)] 
    && \text{(by \ref{slambda en 1-x-y})}\\
    &= \sum_{\lambda\vdash c} \frac{p_\lambda[(1+xy)(1-x)(1-y)]}{z_\lambda} \\
    &= \sum_{\lambda\vdash c} \left(\prod_i\left( 1+(xy)^{\lambda_i}\right)\right)
    \frac{p_\lambda[(1-x)(1-y)]}{z_\lambda} \\
    &= \sum_{k\ge0} (xy)^k\sum_{\lambda\vdash c} \chi^k(\lambda)\frac{p_\lambda}{z_\lambda}[(1-x)(1-y)],
\end{align*}
where 
\[
    \chi^k(\lambda) =
    \begin{cases}
        \sum_{\mu\vdash k}\binom{m_1(\lambda)}{m_1(\mu)}\cdots\binom{m_k(\lambda)}{m_k(\mu)} & \text{for } 0\le k\le \left\lfloor \frac{c}{2}\right\rfloor,\\
        \chi^{c-k}(\lambda) & \text{for } \left\lfloor \frac{c}{2}\right\rfloor< k\le c,\\
        0 &\text{for } c<k.
    \end{cases}
\]

\begin{lem}
We have $\sum\limits_{\lambda\vdash c} \chi^k(\lambda)\frac{p_\lambda}{z_\lambda} = h_{c-k,k}$ for $0\le k\le \left\lfloor\frac{c}{2}\right\rfloor$.
\end{lem}
\begin{proof}
We compute
\begin{align*}
    h_{c-k, k}
    &= s_{c-k}s_k = \sum_{\mu\vdash c-k}\frac{p_{\mu}}{z_\mu}\sum_{\substack{\nu\vdash k}} \frac{p_{\nu}}{z_\nu} = \sum_{\substack{\mu\vdash c-k\\\nu\vdash k}} \frac{p_{\mu\cup\nu}}{z_\mu z_\nu}\\
    &= \sum_{\substack{\lambda\vdash c\\\mu\vdash k}}\frac{p_\lambda}{z_\lambda}\frac{\prod m_i(\lambda)!}{\prod m_i(\mu)! (m_i(\lambda)-m_i(\mu))!} = \sum_{\lambda\vdash c}\chi^k(\lambda)\frac{p_\lambda}{z_\lambda}.
\end{align*}
\end{proof}

Again by Note \ref{note:positive and negative}, we get the following lemma.
\begin{lem}\label{sa en 1-x-y+xy}
For $n\ge2$, the have the following equality: 
\[
s_n[1-x-y+xy] =(1-x)(1-y)\frac{1-(xy)^{n-1}}{1-xy}.
\]
\end{lem}
\begin{proof}
    A row tableau of size $n$ can have at most one entry equal to $-1$, and at most one entry equal to $-2$. The remaining cells $n-2$ can have and any number $0\le k\le n$ of entries equal to 1 and $n-k-2$ entries equal to 2.
\end{proof}
\begin{proof}[Proof of Thm.~\ref{sc o s2 o sa}]
With these lemmas, expression \eqref{eq:1} now becomes

\begin{equation}\label{eq:3}
    \left.\frac{\sum_{k\ge0}^{c} (xy)^k h_{c-k}h_{k}[(1-x)(1-y)]}{(1-x)(1-y)}\right|_{(xy)}
\end{equation}
\begin{equation}\label{eq:4}
    \stackrel{\ref{sa en 1-x-y+xy}}{=} \sum_{k\ge0}^{2c} (xy)^k +\left.\frac{(1-x)(1-y)}{(1-xy)^2} \sum_{k\ge1}^{c-1} (xy)^k-(xy)^{c}-(xy)^{2k}+(xy)^{c+k}\right|_{(xy)}
\end{equation}
\begin{gather*}
    = \frac{1-(xy)^{2c+1}}{1-xy} + \frac{1+xy}{(1-xy)^2}\left( \frac{(xy)-(xy)^{c}}{1-xy} - (n-1)(xy)^{c} \right.
    \hspace{7em}
    \\ \hspace{14em} 
    \left.    -\frac{(xy)^2-(xy)^{2c}}{1-(xy)^2}+\frac{(xy)^{c}-(xy)^{2c}}{1-xy}\right).
\end{gather*}
Let us consider only the terms which affect the first $c$ entries of the hook+column sequence. We obtain
\[
    \frac{1}{1-xy} + \frac{1+xy}{(1-xy)^2}\left(\frac{xy}{1-xy} - \frac{(xy)^2}{1-(xy)^2}\right)
    = \frac{1}{1-xy} + \frac{xy}{(1-xy)^3},
\]
which is precisely the generating function for OEIS A000124, as announced.
So far, we have shown that $s_c\circ s_2\circ s_2$ yields, for $\gamma=0$, a hook+column sequence starting with $(1, 2, 4, \ldots, T_c+1)$, where $T_c$ is the $c$th triangular number.

Furthermore, as it is apparent from \eqref{eq:3}, the coefficient of $(xy)^k$ coincides with the coefficient of $(xy)^{c-k}$, proving that, in fact,
the hook+column sequence of $s_c\circ s_2\circ s_2$ for $\gamma=0$ is equal to \((1, 2, 4, \ldots, T_c+1, T_c+1, \ldots, 4, 2, 1, 0, \ldots)\). In other words,
$ \hcZ{s_c\circ s_2\circ s_a}$ is equal to
\[
 \sum\limits_{k=0}^{2c-1}\min\left\{\frac{k^2+k+2}{2},\frac{(2c-1-k)^2+(2c-1-k)+2}{2}\right\}\cdot s_{(2ac-2k, 2^k)}.
\]

It remains to show that $\hcN{s_c\circ s_2\circ s_a} $ 
is equal to
\[
\sum\limits_{k=1}^{2c-3}\min\left\{\left\lfloor\frac{(k+1)^2}{4}\right\rfloor,\left\lfloor\frac{(2c-1-k)^2}{4}\right\rfloor\right\}\cdot s_{(2ac-2k-1,\ 2^k,\ 1)}.
\]
We only sketch the proof, since the computations are similar. We begin by considering the expression \eqref{eq:4}. But instead of restricting to monomials which are symmetric in $x$ and $y$, we consider the remaining monomials. We obtain this way the generating function for the coefficients in for $\Sigma(s_c\circ s_2\circ s_a,1)$.

This time, the starting sequence is an offset of OEIS A002620 \cite{oeis}, whose general term is $\big\lfloor\frac{(n+1)^2}{4}\big\rfloor$ and which has a generating function $\frac{z}{(1+z)(1-z)^3}$ \cite{plouffe}. This is shown similarly.

This completes the proof of Theorem \ref{sc o s2 o sa}.
\end{proof}

\subsection{Symmetry of hook+column sequences}\label{sec:symmetry}

This section is dedicated to the proof of our main result (Theorem \ref{Symmetry}). Express $s_2$ as $\frac{1}{2}(p_{1,1}+p_2)$ and $s_{1,1}$ as $\frac{1}{2}(p_{1,1}-p_2)$. Using the results and notations of Section \ref{sec:handful}, for each hook+column $\lambda$ we can write 
\begin{gather*}
    2[\lambda]\, (s_2\circ  f)=\#\D{1,1}{} + \sgn_2(\lambda)\#\D{2}{},\\
    2[\lambda] \, (s_{1,1}\circ  f)=\#\D{1,1}{} - \sgn_2(\lambda)\#\D{2}{}.
\end{gather*}

Our aim is to show that if $f$ is flip-symmetric with offset $r$, then $[\lambda] \, (s_2\circ f)$ equals $[\flip(2r-2; \lambda)]\, (s_2\circ f)$. Let us denote $R=2r-2$ and let $\lambda^R := \flip(R; \lambda)$. Note from Table \ref{tab:sign} that the sign function is invariant under the flip involution.
Hence, proving $\#\D{2}{} = \#\D[\lambda^R]{2}{}$ and $\#\D{1,1}{} = \#\D[\lambda^R]{1,1}{}$ will suffice to show the theorem.\\

We begin with well-definedness. Note that for any $r\ge2$, we have $R\ge2$. Fix $\lambda = (\lambda_1, 2^\beta, 1^\gamma)\in\supp(s_2\circ f)$. We claim that the $R$-flip is well defined on $\lambda$. That is, we can find $\delta\ge0$ such that $\lambda_1 = R+2\delta+\gamma$. Indeed, the parity of $2n = \lambda_1+2\beta+\gamma$ implies $\lambda_1-\gamma\equiv R\equiv0 \mod 2$. Moreover, suppose $\mu\in\D{2}{}$, $\lambda_1$ is even and $2\mu_1=\lambda_1$. Then, $r\le\mu_1-m_1(\mu)$ and $\beta\ge2m_2(\mu)$ imply $R\le\lambda_1-\gamma$, as desired. In any other case (e.g., $\mu\in\D{1,1}{0}$), similar computations yield similar results.

The following four facts complete the proof of the theorem:
\begin{itemize}
    \item If $\lambda_1$ is odd, then $\#\D{2}{}=\#\D[\lambda^R]{2}{}$ by Lemma \ref{A=A^r}.
    \item If $\lambda_1$ is even, then $\#\D{2}{0}=\#\D[\lambda^R]{2}{2}$ by Lemma \ref{A.1=A^r.2}. Relabelling $\lambda$ for $\lambda^R$ gives $\#\D[\lambda^R]{2}{0}=\#\D{2}{2}$.
    \item By Lemma \ref{B.1=B^r.2}, $\#\D{1,1}{0}=\#\D[\lambda^R]{1,1}{2}$, and relabelling $\lambda$ for $\lambda^R$ gives $\#\D[\lambda^R]{1,1}{0}=\#\D{1,1}{2}$.
    \item Finally, $\#\D{1,1}{1}=\#\D[\lambda^R]{1,1}{1}$ by Lemma \ref{B.3=B^r.3}.
\end{itemize}

The non-negativity of the resulting coefficients holds from Lemma \ref{s2 o f} and from the fact that Schur positivity is preserved under plethysm by $s_2$.

\begin{lem}\label{A=A^r}
Under the hypotheses of this section, if $\lambda_1$ is odd then $\#\D{2}{}=\#\D[\lambda^R]{2}{}$.
\end{lem}
\begin{proof} If $\lambda_1=R+2\delta+\gamma$ is odd, then so is $\lambda_1^R = R+2\beta+\gamma$. The multiset $\D{2}{}$ only has one element, $\mu$. We have
\[\mu = \left(r+\delta+\frac{\gamma-1}{2},\ 2^{\frac{\beta}{2}},\ 1^{\frac{\gamma-1}{2}}\right)
\xmapsto{r-\text{flip}}
\left(r+\beta+\frac{\gamma-1}{2},\ 2^{\frac{\delta}{2}},\ 1^{\frac{\gamma-1}{2}}\right) =: \mu^r.\]
Note that $\mu^r$ is the only partition appearing in $\D[\lambda^R]{2}{}$. By hypothesis, these two partitions appear with the same multiplicity 
$[\mu]\, f = [\flip(r; \mu)] \, f$ in their corresponding multisets. This completes the proof.
\end{proof} 

The remaining lemmas will follow in the same spirit as the previous one.

\begin{lem}\label{A.1=A^r.2}
Under the hypotheses of this section, if $\lambda_1$ is even then $\#\D{2}{0}=\#\D[\lambda^R]{2}{2}$.
\end{lem}
\begin{proof}
Explicitly, we can write $\D{2}{0}$ and $\D[\lambda^R]{2}{2}$ (up to multiplicities) as follows,
\begin{gather*} 
    \left\{ \left( r + \delta + \frac{\gamma}{2} - 1, \ 2^{m_2},\ 1^{\beta+\frac{\gamma}{2}-2m_2}\right)\ :\ m_2\le\frac{\beta}{2}\right\}
    \\ \text{and}\
    \left\{ \left( r + \beta + \frac{\gamma}{2}, \ 2^{m'_2},\ 1^{\bullet}\right)\ :\ m_2'\le\frac{\delta-1}{2}\right\}.
\end{gather*}
Applying the $r$-flip on every element of $\D{2}{0}$ yields
\[
\flip(r;\D{2}{0}) = \left\{ \left( r + \beta + \frac{\gamma}{2}, \ 2^{m_2'},\ 1^{\beta+\frac{\gamma}{2}-2m_2}\right)\ :\ m_2\le\frac{\beta}{2}\right\}, \]
where \(m_2'=\frac{1}{2}(\delta+2m_2-\beta-1)\).

We aim to identify $\flip(r;\D{2}{0})$ with $\D[\lambda^R]{2}2$. The only thing that remains to show is that $m_2\le\frac{\beta}{2}$ if and only if $m_2'\le\frac{\delta-1}{2}$.
From the expression of $m_2'$, the inequality $m_2'\le\frac{\delta-1}{2}$ simplifies to
$\frac{{\delta-1}+2m_2-\beta}{2}\le{\frac{\delta-1}{2}},$
and it is now clear that this is equivalent to the inequality $m_2\le\frac{\beta}{2}$.
\end{proof} 

\begin{lem}\label{B.1=B^r.2}
Under the hypotheses of this section, we get $\#\D{1,1}{0}=\#\D[\lambda^R]{1,1}{2}$.
\end{lem}
\begin{proof}
Let $\mu^r$ and $\nu^r$ be the $r$-flip of some $(\mu,\nu)\in\D{1,1}{0}$. Note that $|\mu|+|\nu|=|\lambda|$ and $\mu_1+\nu_1=\lambda_1$ imply
\begin{equation}\label{eq:aux}
2m_2(\mu) + m_1(\mu) +  2m_2(\nu) + m_1(\nu) = 2\beta +\gamma.    
\end{equation}
Adding $R=2r-2$ to both sides of the equation, we conclude $\mu_1+\nu_1 = \lambda_1$ if and only if $\mu^r_1+\nu^r_1 = \lambda^R_1+2$.

Again from Equation \eqref{eq:aux}, we get $m_2\le\beta$ if and only if $m_1(\mu)+m_1(\nu)\ge\gamma$.
Knowing that $\mu_1+\nu_1=\lambda_1$, we write
\[\big(r+2m_2(\mu^r)+m_1(\mu)\big)+\big(r+2m_2(\nu^r)+m_1(\nu)\big) = R + 2\delta + \gamma = 2r-2+2\delta+\gamma.
\]
And so, $m_2\le\beta$ if and only if $m_2' + 1 \le \delta$, where $m_2':=m_2(\mu^r)+m_2(\nu^r)$.
In a similar fashion, one can show $\beta\le m_2+m_1$ if and only if $\delta\le m_2'+m_1+1.$

Summing up, we have proved that a pair $(\mu,\nu)$ is in $\D{1,1}{0}$ if and only if $(\mu^r,\nu^r)$ is in $\D[\lambda^R]{1,1}{2}$, and thus $\#\D{1,1}{0}=\#\D[\lambda^R]{1,1}{2}$ by the flip-symmetry hypothesis on $f$.
\end{proof}

\begin{lem}\label{B.3=B^r.3}Under the hypotheses of this section, we get $\#\D{1,1}{1}=\#\D[\lambda^R]{1,1}{1}$.
\end{lem}
\begin{proof} The proof is similar to that of Lemma \ref{B.1=B^r.2}. Following in each step the case in which the equalities are attained, we also get $\chi_{(\mu,\nu)}^\beta = \chi_{(\mu^r,\nu^r)}^\delta$.
\end{proof}

\section{Final comments}
\label{final_comments}

 Theorems \ref{s2 o sb o sa}, \ref{sc o s2 o sa},  and \ref{Symmetry} imply that the iterated plethysm $f = s_2\circ s_2\circ \ldots\circ s_2\circ s_c\circ s_b\circ s_a$ is flip-symmetric when either $b$ or $c$ is equal to 2.  Therefore, the hook+column sequence $\Sigma(f,\gamma) $ is symmetric for each non-negative integer $\gamma$. However, sequences like this one appear to have stronger properties.  Based on our data, we put forth some questions regarding the structural behaviour of hook+column sequences arising from iterated plethysms.

\begin{enumerate}
\item Our first question is a very natural one. 
Are the hook+column sequences of the form $\Sigma(s_{n_1}\circ s_{n_2}\circ\cdots\circ s_{n_k},\gamma)$ symmetric for all $n_1, \ldots, n_k$ and all  $\gamma$?

\item

Our algebraic definition of the flip involution makes sense even if the offset is $r = 0$ or $1$. Our data suggests that if $f$ is flip symmetric with offset $0$ or $1$, then $\Sigma(s_2\circ f,\gamma)$ and $\Sigma(s_{1,1}\circ f,\gamma)$ are also symmetric sequences. But $s_2\circ f$ and $s_{1,1}\circ f$ are not flip-symmetric sequences according to our definition. See for instance Example \ref{ex:counterex}. Based on our analysis of the data, we infer that there is a wider partial symmetry that has yet to be revealed. A more general description of  the flip-symmetry and the flip involution is expected to exist.

\item
A finite sequence $(a_i)_{i=1, \ldots, n}$ is said to be unimodal  if  there exists a $k$ such that  
\[a_0\le a_1\le \cdots \le a_{k-1} \le a_k \ge a_{k+1} \ge \cdots \ge a_{n-1}\ge a_n.
\]
Based on the available data, we ask:  Are the hook+column sequences $\Sigma(s_{n_1}\circ s_{n_2}\circ\cdots\circ s_{n_k},\gamma)$ {unimodal} for all $n_1, \ldots, n_k$ and all $\gamma$?

Another way of rephrasing this question is reminiscent of a celebrated result, the unimodality of the $q$-binomial coefficients. This result can be reduced to deciding whether the following polynomial is unimodal:
\[(s_a\circ s_b)[1+q] = s_a[(b+1)_q] = \binom{a+b}{a}_q=\sum_m b_m q^m,\]
where $(n)_q$ is the $q$-analogue of $n$ and $\binom{n}{k}_q$ is the $q$-binomial coefficient. On the other hand, our problem can be reduced to studying whether
\[\frac{s_{n_1}\circ s_{n_2}\circ\cdots\circ s_{n_k}[1-x-y]\cdot (x-y)}{(-1)^\gamma(1-x)(1-y)(x^{\gamma+1}y^{\gamma+1})} = \sum_\beta a_\beta(xy)^\beta\]
is a unimodal polynomial in the variable $(xy)$ for any $\gamma$.

\item
A positive finite sequence $(a_i)_{i=1, \ldots, n}$ is said to be log-concave if $a_i^2 \ge a_{i-1}a_{i+1}$ for all $i = 2, \ldots,  n-1$. It is well known that log-concavity implies unimodality. Recently, the log-concavity of many combinatorial sequences has been established thanks to the development of several breakthrough methods \cite{branden14,branden19,baker17,braden20}. In some of these works, the sequences shown to be log-concave are not the classical combinatorial sequences, but a renormalization of them \cite{huh22}.

The hook+column sequence $\Sigma(s_2^{\circ5},\gamma)$ appearing in Figure \ref{fig:normalidad} gives us an example of a 
hook+column sequence that is not log-concave.
However, is there a sensible renormalization of hook+column sequences arising from plethysm that renders them log-concave?

\item
Asymptotic normality is another structural phenomenon commonly found in combinatorial sequences. Experimental evidence suggests that the hook+column sequences of $(s_2)^{\circ k}:=s_2\circ \stackrel{k\text{ times}}{\ldots}\circ s_2$ for any fixed $\gamma$ are asymptotically normal when $k$ tends to infinity, as Figure \ref{fig:normalidad} illustrates. Moreover, a $\chi^2$ normality test returns a $p$-value of 1 or almost 1 for every sequence coming from $s_2^{\circ k}$, $k=2,\ldots,5$. These huge $p$-values\footnote{We adopt the usual convention of rejecting the null hypothesis if the $p$-value is smaller than 0.05.} seem to indicate that the Gaussian curve \emph{perfectly} fits our sequences, even for small values of $k$.

\begin{figure}[H]\label{table_normal}
    \centering
    \footnotesize 
    \begin{tabular}{  l | l | l }
                \toprule
                Function, $f$ & $\gamma$  & Hook+column sequence, $\Sigma(f, \gamma)$  \\ 
                \midrule
                $s_2^{\circ2}$ & 0 & (1, 1)\\
                $s_2^{\circ3}$ & 0 & (1, 2, 2, 1)\\
                $s_2^{\circ4}$ & 0 & (1, 3, 8, 13, 13, 8, 3, 1)\\
                $s_2^{\circ5}$ & 0 & (1, 4, 20, 72, 205, 446, 756, 986, 986, 756, 446, 205, 72, 20, 4, 1)\\
                \bottomrule
            \end{tabular}
        \includegraphics[width=0.8\textwidth]{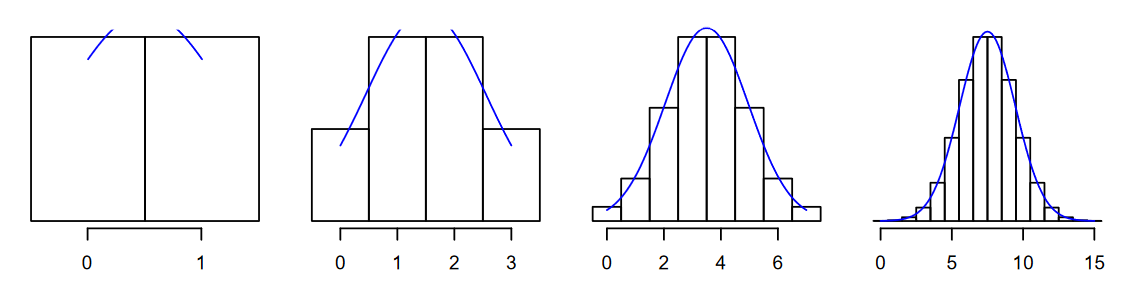}
    \caption{On top, a table showing the hook+column sequence of $s_2^{\circ k}$ for $\gamma=0$ and $k=2,3,4,5$. Below, plots of the aforementioned sequences with the $x$-axis being $\beta$, and represented as the normalized histogram whose frequencies read $\Sigma(s_2^{\circ k}, 0)$. They appear overlaid with Gaussian curves of adjusted mean and variance.}
    \label{fig:normalidad}
\end{figure}

Is the hook+column sequence $\Sigma({s_2^{\circ k}},\gamma)$ asymptotically normal for each fixed $\gamma$? (In the sense that its relative sums approach a Gaussian curve when $k$ tends to infinity.)
(See \cite{billey, charalambides, harper} for more details in asymptotic normality of combinatorial integer sequences.)

\item
We have used SageMath \cite{sage} to compute   data supporting these questions. For instance, the hook+column sequences of the family $s_{1,1}^{\circ k}$ also appear to be asymptotically normal. Moreover, our data for $f_{abc} := s_c\circ s_b\circ s_a$ suggest that the limiting hook+column sequences of $f_{abc}$ when both $b$ and $c$ tend to infinity is asymptotically normal (see Figure \ref{fig:normalPlots}).

\begin{figure}[H]
    \centering
    \includegraphics[width=0.8\textwidth]{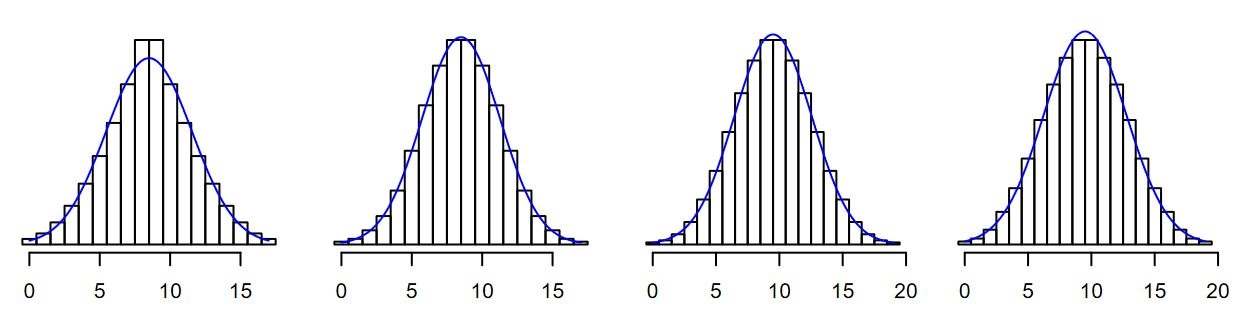}
    \caption{From left to right, the histogram plots for the hook+column sequences associated to $s_9\circ s_2\circ s_2$, $s_6\circ s_3\circ s_2$, $s_5\circ s_4\circ s_2$ and $s_4\circ s_5\circ s_2$, with $\gamma=0$, and where the $x$-axis represents $\beta$. They appear overlaid with Gaussian curves of adjusted mean and variance.}
    \label{fig:normalPlots}
\end{figure}

 \item In Theorems \ref{s2 o sb o sa} and \ref{sc o s2 o sa}, we gave explicit descriptions of the hook+column sequences $\Sigma(f_{abc}, \gamma)$, whenever $b$ or $c$ are equal to 2. Our data allows us to make reasonable guesses about what the hook+column sequences of $f_{abc}$ approach for other values of $b$ and $c$. (Recall that, by Lemma \ref{tail0}, the value of $a$ does not affect the non-vanishing part of the sequence.)
 
The hook+column sequences of $f_{23c}=s_c\circ s_3 \circ s_2$ and $\gamma=0$ up to $c=6$ are shown in Table \ref{tab:snos3os2}. 
Unlike the sequences in Theorems \ref{s2 o sb o sa}
and \ref{sc o s2 o sa}, each consecutive sequence is not simply a longer version of the previous ones. However, they seem to stabilize. Is their stable limit sequence (1, 2, 5, 10, 19, 33, 57, 92, 147, \ldots),
the number of partitions with two kinds of 1s, 2s, and 3s? (OEIS A000098 \cite{oeis}.)

The hook+column sequences of $f_{24c}=s_c\circ s_4 \circ s_2$ and $\gamma=0$ up to $c=6$ are shown in Table \ref{tab:snos4os2}. Again, the coefficients seem to stabilize. Is their stable limit sequence (1, 2, 5, 11, 22, 42, 77, 135, \ldots),
the number of partitions of $2n$? (OEIS A058696 \cite{oeis}.)

More generally, are the hook+column sequences $\Sigma(f_{abc}, \gamma)$ counting partitions of a number with some restriction  on the
allowable parts, for all $b, c$ and $\gamma$? 
The literature already contains instances of sequences similar to the last two, in the context of exploring structural constants of symmetric functions, including Kronecker coefficients \cite{colmenarejo15}.
\end{enumerate}

\section*{Acknowledgements}

The authors express their appreciation to Adrià Lillo, Emmanuel Briand, and Laura Colmenarejo for their insightful comments and engaging discussions. They also acknowledge the outstanding efforts of the two anonymous referees in thoroughly reviewing the preliminary version of their work.

\bibliography{BibliographyPartial}

\end{document}